\documentclass{amsart}
\usepackage{amssymb, amsbsy, amsthm, amsmath, amstext, amsopn}
\usepackage[all]{xy}
\usepackage{amsfonts}
\usepackage{amscd}

\newtheorem{thm}{Theorem} 

\newtheorem{prop}[thm]{Proposition}

\theoremstyle{definition}

\newtheorem{defn}[thm]{Definition}

\theoremstyle{remark}

\begin{document}

\newcommand{\hs}{\mbox{\hspace{.4em}}}
\newcommand{\ds}{\displaystyle}
\newcommand{\bd}{\begin{displaymath}}
\newcommand{\ed}{\end{displaymath}}
\newcommand{\bcd}{\begin{CD}}
\newcommand{\ecd}{\end{CD}}

\newcommand{\proj}{\operatorname{Proj}}
\newcommand{\bproj}{\underline{\operatorname{Proj}}}
\newcommand{\spec}{\operatorname{Spec}}
\newcommand{\bspec}{\underline{\operatorname{Spec}}}
\newcommand{\pline}{{\mathbf P} ^1}
\newcommand{\pplane}{{\mathbf P}^2}

\newcommand{\ldb}{[[}
\newcommand{\rdb}{]]}

\newcommand{\Sym}{\operatorname{Sym}^{\bullet}}
\newcommand{\Symp}{\operatorname{Sym}}

\newcommand{\cA}{{\mathcal A}}
\newcommand{\bA}{{\mathbf A}}
\newcommand{\cB}{{\mathcal B}}
\newcommand{\cC}{{\mathcal C}}
\newcommand{\cs}{{\mathbf C} ^*}
\newcommand{\boldc}{{\mathbf C}}
\newcommand{\cE}{{\mathcal E}}
\newcommand{\cF}{{\mathcal F}}
\newcommand{\cK}{{\mathcal K}}
\newcommand{\cL}{{\mathcal L}}
\newcommand{\M}{{\mathcal M}}
\newcommand{\bM}{{\mathbf M}}
\newcommand{\theo}{\mathcal{O}}
\newcommand{\boldp}{{\mathbf P}}
\newcommand{\boldq}{{\mathbf Q}}
\newcommand{\cQ}{{\mathcal Q}}
\newcommand{\cU}{{\mathcal U}}

\newcommand{\End}{\operatorname{End}}
\newcommand{\Hom}{\operatorname{Hom}}
\newcommand{\uHom}{\underline{\operatorname{Hom}}}
\newcommand{\Ext}{\operatorname{Ext}}
\newcommand{\bExt}{\operatorname{\bf{Ext}}}
\newcommand{\uTor}{\underline{\operatorname{Tor}}}

\newcommand{\inv}{^{-1}}
\newcommand{\airtilde}{\widetilde{\hspace{.5em}}}
\newcommand{\airhat}{\widehat{\hspace{.5em}}}
\newcommand{\nt}{^{\circ}}

\newcommand{\supp}{\operatorname{supp}}
\newcommand{\id}{\operatorname{id}}
\newcommand{\res}{\operatorname{res}}
\newcommand{\lrar}{\leadsto}
\newcommand{\im}{\operatorname{Im}}

\newcommand{\TF}{\operatorname{TF}}
\newcommand{\Bun}{\operatorname{Bun}}
\newcommand{\Hilb}{\operatorname{Hilb}}
\newcommand{\nthord}{^{(n)}}
\newcommand{\Aut}{\underline{\operatorname{Aut}}}
\newcommand{\Gr}{\operatorname{\bf Gr}}
\newcommand{\Fr}{\operatorname{Fr}}

\title[moduli of framed sheaves]{Representability for some moduli stacks of framed sheaves}
\author{Thomas A. Nevins}
\address{Department of Mathematics\\University of Michigan\\Ann Arbor, MI 48109-1109 USA}
\email{nevins@umich.edu}
\subjclass{Primary 14D22; Secondary 14D20}


\maketitle

\noindent
{\bf Introduction.} \;\;
Moduli problems for various kinds of framed sheaves have been studied and used in many settings (see, for example, \cite{MR95e:14006}, \cite{MR91m:32031}, \cite{MR95i:53051}), and there is a good general theory of moduli for semistable framed sheaves, thanks to the work of Huybrechts and Lehn (\cite{MR96a:14017}, \cite{MR95i:14015}).  By contrast, there seem to be only a few examples in which the {\em full} moduli functor for framed sheaves (without conditions of semistability) is known to be represented by a scheme.  In this paper, we prove a representability theorem for the full moduli functors of framed torsion-free sheaves on projective surfaces under certain conditions.

\vspace{.5em}

Let $S$ denote a smooth, connected complex projective surface, and let \mbox{$D\subset S$} denote a smooth connected complete curve in $S$.  Fix a vector bundle $E$ on $D$.  An {\em $E$-framed torsion-free sheaf on $S$} is a pair $(\cE,\phi)$ consisting of a torsion-free sheaf $\cE$ on $S$ and an isomorphism $\phi: \cE\big|_{D}\rightarrow E$; the isomorphism $\phi$ is called an {\em $E$-framing} of $\cE$.  An {\em isomorphism} of $E$-framed torsion-free sheaves on $S$ is an isomorphism of the underlying torsion-free sheaves on $S$ that is compatible with the framings.  Let $\TF_S(E)$ denote the moduli functor for isomorphism classes of $E$-framed torsion-free sheaves on $S$.  The reader should note that in the work of Huybrechts--Lehn the framing $\phi$ need {\em not} be an isomorphism; as a consequence of our more restrictive definition, the moduli functors that we study have no hope of being proper.

\vspace{.5em}

Suppose the vector bundle $E$ satisfies
\begin{equation}\label{poscondition}
H^0\left(D, \End E\otimes N^{-k}_{D/S} \right) = 0
\end{equation}
for all $k\geq 1$; here $N_{D/S}$ is the normal bundle of $D$ in $S$.  If $D\subset S$ is an arbitrary curve, there may be very few such bundles.  However, if $D$ is smooth and has positive self-intersection in $S$, then $N_{D/S}\inv$ is a negative line bundle on $D$, and consequently this condition on $E$ is an open condition which is satisfied by all semistable vector bundles on $D$.
\begin{thm}\label{maintheorem}
Suppose that $S$ is a smooth, connected complex projective surface and $D\subset S$ is a smooth connected complete curve. Suppose, in addition, that $E$ is a vector bundle on $D$ that satisfies Condition \eqref{poscondition} for all $k\geq 1$.  Then the functor $\TF_S(E)$ is represented by a scheme. 
\end{thm}
\noindent
In the proof of Theorem \ref{maintheorem} we work in the slightly more general setting of a family of vector bundles on $D$, parametrized by a scheme $U$, that satisfies Condition \eqref{poscondition} for all $k\geq 1$ at every point of $U$.  Note also that the reader who is familiar with the language of stacks may restate Theorem \ref{maintheorem} in the following form: over the substack of its target that parametrizes vector bundles on $D$ that satisfy Condition \eqref{poscondition} for all $k\geq 1$, the fibers of the restriction morphism from the moduli stack of torsion-free sheaves on $S$ that are locally free along $D$ to the moduli stack of vector bundles on $D$ are schemes.

\vspace{.5em}

Functors of the type we study here arose naturally (in some special cases) in the representation-theoretic constructions of Nakajima; Theorem \ref{maintheorem} demonstrates that the existence of the fine moduli schemes used by Nakajima is a much more general phenomenon, one which we hope can be exploited more widely in the study of sheaves on noncompact surfaces.
The new ingredient in our proof of Theorem \ref{maintheorem} is the use of formal geometry along the curve $D$; in particular, the techniques used here are completely different from those of \cite{MR96a:14017}, \cite{MR95i:14015}, and make no use of geometric invariant theory (GIT).  Although Lehn (\cite{MR95g:14013}) has, under some conditions on the curve $D$ and the bundle $E$ along the curve, proven that the full moduli functors for vector bundles on $S$ with framing along $D$ by $E$ are represented by {\em algebraic spaces},  from the point of view of the usual GIT techniques it is perhaps surprising that there is a fine moduli {\em scheme} (a much stronger fact) for all framed sheaves: indeed, there can be framed sheaves that are not semistable for {\em any} polarization.

\vspace{.5em}

The author is grateful to K. Corlette for his constant guidance and also for some essential conversations. The author is grateful also to {V. Baranovsky}, T. Pantev and I. Robertson for valuable discussions related to this work, and to a referee for helpful suggestions (especially for improvements in the proof of Proposition \ref{helped}).  The author's graduate work at the University of Chicago, of which this paper is a result, was supported in part by an NDSEG fellowship from the Office of Naval Research.

\vspace{.8em}

\noindent
{\bf Affine bundles over $\Bun\nt(D)$.}\label{affinebundlessection}\;\;
In this section we construct the fundamental affine bundles $\bA_n$ (for $n$ in the range $1\leq n < \infty$) over $U$ that we will use to embed
the functor $\TF_S(E)$ in a scheme.  The construction of these bundles and the
description of the universal properties they possess must be well known (cf. \cite{MR34:6796}, in which the relevant cohomology groups are discussed), but
the author does not know a suitable reference.  

Fix a surface $S$, a curve $D$ in $S$, a scheme $U$, and a vector bundle $E$ on $D\times U$ as in Theorem \ref{maintheorem}.  Let $D\nthord$ (that is, $D$ with structure sheaf $\theo_S/I_D^{n+1}$,
 $0\leq n <\infty$) denote the $n$th order neighborhood of $D$ in $S$.

\begin{defn}
Let $\cA_n$ denote the moduli functor over $U$ of isomorphism classes of triples $(\cE,V\xrightarrow{f} U, \phi)$ consisting of 
\begin{enumerate}
\item a vector bundle $\cE$ on $D\nthord\times V$,
\item a morphism $f:V\rightarrow U$, and
\item an isomorphism $\phi: \cE\big|_{D\times V} \rightarrow (1_D\times f)^*E.$
\end{enumerate}
\end{defn}
\noindent
Suppose that $\cE$ is a vector bundle over $D\nthord$; then $\cE$ has a canonical (decreasing) filtration as an $\theo_{D\nthord}$-module with filtered pieces
$\ds F_j\cE = I_D^j\cE$,
where $I_D$ is the ideal of $D\subset D\nthord$.
By its construction, this filtration is preserved by any endomorphism of the vector bundle $\cE$, and moreover
$F_j\cE/F_{j+1}\cE \cong N_{D/S}^{-j} \otimes \big(F_0\cE/F_1\cE\big)$ provided $0\leq j\leq n$.
Using these facts together with the exact sequence
\bd
0 \rightarrow \Hom (E, E\otimes N_{D/S}^{-n}) \rightarrow \Hom(\cE,\cE) \rightarrow \Hom \big(\cE, \cE\big|_{D^{(n-1)}}\big) \rightarrow 0
\ed
and condition \eqref{poscondition},
one may prove by induction on $n$ that $\End(\cE) \subseteq \End(\cE\big|_D)$ and consequently that $E$-framed bundles on $D\nthord$ are rigid.

Evidently $\cA_0 \cong U$; moreover, there are maps
$\ds\pi_{n+1}: \cA_{n+1} \rightarrow \cA_n$
for all $n\geq 0$.

\begin{prop}\label{helped}
Each $\cA_n$  ($n\geq 1$) is represented by a scheme $\bA_n$ that is an 
affine bundle over $\bA_{n-1}$.
\end{prop}

\begin{proof}
Working inductively, it will suffice to construct an $\bA_{n-1}$-scheme $\bA_n$ that represents $\cA_n$ and is an affine bundle over $\bA_{n-1}$.  
Fix a universal bundle $E^{(n-1)}$ on $D^{(n-1)}\times\bA_{n-1}$.  For any scheme $T$, an element of $\cA_n(T)$ determines a map $f: T\rightarrow \bA_{n-1}$, and, if $(\cE,\phi)$ is the given element of $\cA_n(T)$, there is an isomorphism of $\cE\big|_{D^{(n-1)}\times T}$ with $(1\times f)^*E^{(n-1)}$ compatibly with the framings by $E$.  But then, because $E$-framed bundles on $D^{(n-1)}$ are rigid, we find that $\cA_n$ as a functor over $\bA_{n-1}$ is isomorphic to the functor taking $f: T\rightarrow\bA_{n-1}$ to the set of isomorphism classes of pairs $(\cE,\phi)$ consisting of a bundle $\cE$ on $D\nthord\times T$ together with an isomorphism $\phi$ of $\cE\big|_{D^{(n-1)}\times T}$ with $(1\times f)^*E^{(n-1)}$.  We will refer to such a pair as an $E^{(n-1)}$-framed bundle.

Because the statement of the proposition is local on $\bA_{n-1}$, we may assume that $\bA_{n-1}$ is an affine scheme that is the spectrum of a local ring $R$.
For simplicity, write $\theo = \theo_{D^{(n-1)}\times\bA_{n-1}}$ and $\theo' = \theo_{D^{(n)}\times\bA_{n-1}}$.  The ``change of rings'' spectral sequence (see Chap. XVI, Section 5 of \cite{MR17:1040e})
\bd
E_2^{p,q} = \Ext_{\theo}^p\big(\uTor_q^{\theo'}(E^{(n-1)},\theo), E(-nD)\big) \Rightarrow \Ext_{\theo'}^{p+q}\big(E^{(n-1)}, E(-nD)\big)
\ed
yields the exact sequence of terms of low degree
\begin{multline}\label{termsoflowdegree}
0 \rightarrow \Ext^1_{\theo}\big(E^{(n-1)}, E(-nD)\big) \rightarrow \Ext^1_{\theo'}\big(E^{(n-1)}, E(-nD)\big)\\
 \xrightarrow{\beta} \Hom\big(\uTor_1^{\theo'}(E^{(n-1)}, \theo), E(-nD)\big) \rightarrow 0.
\end{multline}
Note that $\beta$ is surjective since the next term in the sequence is $\Ext^2_{\theo}\big(E^{(n-1)}, E(-nD)\big)$, which vanishes because $D$ is one-dimensional.  Using $\uTor_1^{\theo'}(E^{(n-1)},\theo) \cong E(-nD)$ one may check that there is a canonical element $e$ of $\Hom\big(\uTor_1^{\theo'}(E^{(n-1)}, \theo), E(-nD)\big)$ such that $\beta\inv(e)$ is exactly the $\Ext^1_{\theo}\big(E^{(n-1)}, E(-nD)\big)$-subtorsor of $\Ext^1_{\theo'}\big(E^{(n-1)}, E(-nD)\big)$ that classifies 1-extensions 
\bd
0 \rightarrow E(-nD) \rightarrow \cE \rightarrow E^{(n-1)} \rightarrow 0
\ed
for which $\cE$ is a locally free $\theo'$-module.  
Now,
Condition \eqref{poscondition}, together with Cohomology and Base Change, implies that the $R$-module $\Ext^1_{\theo}\big(E^{(n-1)}, E(-nD)\big) \cong H^1\big(D\times \bA_{n-1},\End(E)\otimes N^{-n}_{D/S}\big)$ is projective, hence free.  One can easily construct, moreover, a universal 1-extension over $D^{(n)}\times\bA_{n-1}\times\beta\inv(e)$ (using, for example, an affine subspace of the \v{C}ech cocycles that maps isomorphically to $\beta\inv(e)$ to furnish gluing data).  Because the exact sequence \eqref{termsoflowdegree} and the element $e$ are functorial under pullback along morphisms of affine schemes $\spec R' \xrightarrow{f} \spec R = \bA_{n-1}$, this universal 1-extension induces a functorial bijection between the set $\beta_{R'}^{-1}(e)$ (the inverse image of the canonical element under the base-changed map $\beta$) and the set of isomorphism classes of pairs $(\cE,\phi)$ consisting of a vector bundle $\cE$ on $D^{(n)}\times\spec R'$ and a framing $\phi: \cE\big|_{D^{(n-1)}\times\spec R'} \rightarrow (1\times f)^*E^{(n-1)}$.

 Consequently $\cA_n$ is represented as a functor over $\bA_{n-1}$ by the torsor over $\spec\Sym\Ext^1_{\theo}\big(E^{(n-1)}, E(-nD)\big)$ defined by $\beta\inv(e)$, proving the proposition.
\end{proof}

\vspace{.8em}

\noindent
{\bf Proof of Theorem \ref{maintheorem}.}\label{proofsection}\;\;
There is a compatible family of morphisms 
$F_n: \TF_S(E) \rightarrow \bA_n$
given by restriction.
Fix a $\spec\boldc$-valued point of $\TF_S(E)$, that is, a point $u\in U$ together with an $E_u$-framed pair $(\cF, \phi)$ on $S$.  We will show that there is an open subfunctor $Z$ of $\TF_S(E)$ that contains $(\cF,\phi)$ and is represented by a scheme.

Fix a polarization $H$ of $S$, and choose $m$ sufficiently large that
\begin{enumerate}
\item $\cF\otimes H^m$ is globally generated and 
\item $H^1(\cF\otimes H^m) = H^2 (\cF\otimes H^m) = 0$.
\end{enumerate}
Further, fix $n$ sufficiently large that the restriction map
\bd
H^0\left(\cF\otimes H^m\right)\rightarrow H^0\left(\cF\otimes H^m\big|_{D\nthord}\right)
\ed
is injective; it is possible to choose such an $n$ because $\cF$ is torsion-free.  Finally, choose 
$m'$ sufficiently large that $\ds H^1\left(\cF\otimes H^{m+m'}\big|_{D\nthord}\right) =0$.

Next, let $Z\subseteq \TF_S(E)$ denote the open subfunctor parametrizing those triples \newline $\ds\left(W\xrightarrow{f} U, \cE, \phi: \cE\big|_{D\times W}\rightarrow (1\times f)^*E\right)$ for which the family $\cE$ satisfies the following conditions:
\begin{enumerate}
\item[a.] $\cE_w\otimes H^m$ is globally generated for all $w\in W$,
\item[b.] $H^1(\cE_w\otimes H^m) = H^2(\cE_w\otimes H^m) = 0$ for all $w\in W$,
\item[c.] the map $\ds H^0\left(\cE_w\otimes H^m\right)\rightarrow H^0\left(\cE_w\otimes H^m\big|_{D\nthord}\right)$ is injective for all $w\in W$, and
\item[d.] $H^1\left( \cE_w\otimes H^{m+m'}\big|_{D\nthord}\right) = 0$ for all $w\in W$.
\end{enumerate}

\noindent
In the previous section we showed that there is a universal vector bundle $E\nthord$ on $D\nthord\times \bA_n$.  Fix an element of $Z(W)$; then the map
$F_n(W): W\rightarrow \bA_n$
yields a vector bundle $(1\times F_n)^*E\nthord$ on $D\nthord\times W$ together with an isomorphism
\bd
\cE_W\big|_{D\nthord\times W} \xrightarrow{\phi_n} (1\times F_n)^*E\nthord;
\ed
here $\cE_W$ denotes the torsion-free sheaf on $S\times W$ determined by the fixed element of $Z(W)$.  Let $p_W$ denote the projection $S\times W\rightarrow W$.  Then by construction the sheaves $(p_W)_* \cE_W\otimes H^m$, $(p_W)_*\cE_W\otimes H^{m+m'}$, and $(p_W)_*\left(\cE_W\otimes H^{m+m'}\big|_{D\nthord\times W}\right)$ are vector bundles on $W$, and, choosing a section $s$ of $H^{m'}$ the zero locus of which has transverse intersection with $D$, there is a commutative diagram
\bd
\xymatrix{(p_W)_*\cE_W\otimes H^m \ar[r]\ar[d]^{\otimes s} & (p_W)_*\left( \cE_W\otimes H^m\big|_{D\nthord\times W}\right)\ar[d]^{\otimes s}\\
(p_W)_*\cE_W\otimes H^{m+m'} \ar[r] & (p_W)_*\left( \cE_W\otimes H^{m+m'}\big|_{D\nthord\times W}\right)}
\ed
for which the vertical arrows (given by tensoring with $s$) and the top row are injective.  Using $\phi_n$, we may replace this diagram canonically with the diagram
\bd
\xymatrix{(p_W)_*(\cE_W\otimes H^m)\ar[r] \ar[d]^{\otimes s} & (p_W)_* \left( (1\times F_n)^*E\nthord\otimes H^m\right) \ar[d]^{\otimes s} \\
(p_W)_*\cE_W\otimes H^{m+m'} \ar[r] & (p_W)_*\left( (1\times F_n)^* E\nthord\otimes H^{m+m'}\right).}
\ed

Now, by assumption (d) on $W$, we have
\bd
(p_W)_*\left( (1\times F_n)^*E\nthord \otimes H^{m+m'}\right) = F_n^*\left( (p_{\bA_n})_*(E\nthord\otimes H^{m+m'})\right),
\ed
where $p_{\bA_n}: D\nthord\times\bA_n \rightarrow \bA_n$ is the projection, and so finally we obtain the diagram of vector bundles
\bd
\xymatrix{(p_W)_*\cE_W\otimes H^m \ar[d]^{\otimes s} \ar[dr]^{r} & \\
(p_W)_*\cE_W\otimes H^{m+m'} \ar[r] & F_n^*\left( (p_{\bA_n})_* (E\nthord\otimes H^{m+m'})\right)}
\ed
on $W$, where the diagonal map $r$ and the map $\otimes s$ are injective.  By construction, 
furthermore, the image of the morphism $r$ is a vector subbundle of 
$F_n^*\left( (p_{\bA_n})_*(E\nthord \otimes H^{m+m'})\right)$ and consequently 
determines a morphism $W\rightarrow \Gr$ over $\bA_n$, where
$\Gr \xrightarrow{q} \bA_n$
denotes the relative Grassmannian for the vector bundle $(p_{\bA_n})_*(E\nthord\otimes H^{m+m'})$ on $\bA_n$, the fiber of which over $a\in\bA_n$ parametrizes vector subspaces of $H^0(E\nthord\otimes H^{m+m'})$ that are of dimension $h^0(\cF\otimes H^m)$.

We now construct a Quot-scheme over $\Gr$ that we will use to represent $Z$.  We may pull back $(p_{\bA_n})_*(E\nthord\otimes H^{m+m'})$ to $\Gr$ to obtain a vector bundle $q^*(p_{\bA_n})_*(E\nthord\otimes H^{m+m'})$ on (an open subset of) $\Gr$, with universal subbundle
\bd
\cU \subset q^*(p_{\bA_n})_*(E\nthord\otimes H^{m+m'})
\ed
of rank $h^0(\cF\otimes H^m)$.  If $p_{\Gr}: S\times \Gr \rightarrow \Gr$ denotes the projection to $\Gr$, we obtain a bundle
$\ds p_{\Gr}^*\cU \subset p_{\Gr}^*q^*(p_{\bA_n})_*(E\nthord\otimes H^{m+m'})$
on $S\times \Gr$, as well as a quotient
\bd
p_{\Gr}^*q^*(p_{\bA_n})_* (E\nthord\otimes H^{m+m'})\rightarrow (1\times q)^*(E\nthord\otimes H^{m+m'})
\ed
and subquotient
$\ds (1\times q)^*(E\nthord\otimes H^m) \subset (1\times q)^*(E\nthord\otimes H^{m+m'})$
that are sheaves on $S\times\Gr$ supported on $D\nthord\times\Gr$.

Consider the relative Quot-scheme
$q':\operatorname{Quot}_{S\times\Gr/S}(p_{\Gr}^*\cU) \longrightarrow \Gr$
that parametrizes quotient sheaves for the family $p_{\Gr}^*\cU$ on $S\times\Gr /S$.  There is a universal quotient $(1\times q')^*p_{\Gr}^*\cU \rightarrow\cQ$ on $S\times \operatorname{Quot}_{S\times\Gr/S}$, giving a diagram
\begin{equation}\label{asteriskdiagram}
\begin{split}
\xymatrix{(1\times q')^*p_{\Gr}^*\cU \ar[r]\ar[d] & (1\times q')^*p_{\Gr}^*q^*(p_{\bA_n})_*(E\nthord\otimes H^{m+m'})\ar@<8ex>[d]\\
\cQ & (1\times qq')^*(E\nthord \otimes H^m) \subset (1\times qq')^*(E\nthord\otimes H^{m+m'}).}
\end{split}
\end{equation}
There is a closed subscheme of $\operatorname{Quot}_{S\times\Gr/S}$ (see the proof of Theorem 1.6 of \cite{MR88b:14006}) that represents the subfunctor of those quotients the kernels of which project to zero in 
$(1\times qq')^*(E\nthord\otimes H^{m+m'})$, and a closed subscheme $\cC$ of that closed subscheme that represents the sub-subfunctor that parametrizes those quotients the images of which in $(1\times qq')^*(E\nthord\otimes H^{m+m'})$ actually lie in the subsheaf $(1\times qq')^*(E\nthord\otimes H^m)$.  $\cC$ then represents the functor of quotients of $p_{\Gr}^*\cU$ that map to $(1\times qq')^*(E\nthord\otimes H^m)$---that is, it is exactly the closed subscheme over which Diagram \eqref{asteriskdiagram} extends to
\begin{equation}
\begin{split}
\xymatrix{(1\times q')^*p_{\Gr}^*\cU \ar[r]\ar[d] & (1\times q')^*p_{\Gr}^*q^*(p_{\bA_n})_*(E\nthord\otimes H^{m+m'})\ar@<8ex>[d]\\
\cQ\ar[r] & (1\times qq')^*(E\nthord \otimes H^m) \subset (1\times qq')^*(E\nthord\otimes H^{m+m'}).}
\end{split}
\end{equation}
Restricting further to an open subscheme $\cC\nt$ of $\cC$, we may assume that, over $\cC\nt$, the map
$\ds\cQ\big|_{D\nthord\times \cC\nt} \rightarrow (1\times qq')^*(E\nthord\otimes H^m)$
is an isomorphism, that $\cQ$ is a family of torsion-free sheaves on $S$, and that conditions (a) through (d) are satisfied.  

By construction the morphism $W\rightarrow \Gr$ lifts to a morphism $W\rightarrow \cC\nt$; this construction thus determines a morphism of functors $Z\rightarrow \cC\nt$.  Similarly, there is a forgetful morphism $\cC\nt\rightarrow Z$.  Finally, it is clear from the construction that these two morphisms of functors are inverses of each other, as desired. \qedsymbol

\bibliographystyle{alpha}

\end{document}